# A powerful test based on tapering for use in functional data analysis

**Dan J. Spitzner**

*Department of Statistics*
*University of Virginia*
*P. O. Box 400135*
*Charlottesville, VA 22904-4135*
*USA*
*e-mail:* spitzner@virginia.edu

**Abstract:** A test based on tapering is proposed for use in testing a global linear hypothesis under a functional linear model. The test statistic is constructed as a weighted sum of squared linear combinations of Fourier coefficients, a tapered quadratic form, in which higher Fourier frequencies are down-weighted so as to emphasize the smooth attributes of the model. A formula is $Q_n^{OPT} = n \sum_{j=1}^{p_n} j^{-1/2} \|\boldsymbol{Y}_{n,j}\|^2$. Down-weighting by $j^{-1/2}$ is selected to achieve adaptive optimality among tests based on tapering with respect to its "rates of testing," an asymptotic framework for measuring a test's retention of power in high dimensions under smoothness constraints. Existing tests based on truncation or thresholding are known to have superior asymptotic power in comparison with any test based on tapering; however, it is shown here that high-order effects can be substantial, and that a test based on $Q_n^{OPT}$ exhibits better (non-asymptotic) power against the sort of alternatives that would typically be of concern in functional data analysis applications. The proposed test is developed for use in practice, and demonstrated in an example application.

## Contents

## 1. Introduction

The subject of this article is the functional linear model and functional linear hypothesis, both cornerstones functional data analysis (FDA) methodologies. This model is described as a sample of independent random functions (sometimes called curves or profiles), which here will be taken to have a common one-dimensional domain, $(a, b]$, and real-valued response. To facilitate study of asymptotic properties, the present investigation will adopt a *replicated* functional linear model, in which each of $N$ response points of dimension $P$ is replicated $n$ times. In matrix form, the $i$'th replication is

$$d\boldsymbol{Y}_i(t) = \boldsymbol{X}\boldsymbol{\beta}(t)dt + \sigma d\boldsymbol{\epsilon}_i(t), \quad (1.1)$$

for $i = 1, \ldots, n$, where $\boldsymbol{Y}_i(t) = [Y_{i,1}(t), \ldots, Y_{i,N}(t)]^T$ is a functional vector of responses on $t \in (a, b]$, and $\boldsymbol{X} = [\boldsymbol{x}_1, \ldots, \boldsymbol{x}_N]^T$ is a $N \times P$ "essence" regressor matrix, assumed of full column-rank, which stores the values of explanatory variables at which the response functions are measured. The vector $\boldsymbol{\beta}(t) = [\beta_1(t), \ldots, \beta_P(t)]^T$ is a $P \times 1$ functional vector of regression coefficients, and the $\boldsymbol{\epsilon}_i$ are independent and identically distributed functional error vectors; each $\boldsymbol{\epsilon}_i(t) = [\epsilon_{i,1}(t), \ldots, \epsilon_{i,N}(t)]^T$ is a vector of independent and identically distributed error functions, for which $E[\epsilon_{i,k}(t)] = 0$ and $V[\epsilon_{i,k}(t)] = 1$ for each $(i, k)$ and $t$, and which are stationary on $t \in (a, b]$. The functional linear hypothesis is $H_0 : \boldsymbol{L}^T\boldsymbol{\beta}(t) = 0$ for all $t \in (a, b]$ against a general alternative, where $\boldsymbol{L}$ is a $P \times \nu$ hypothesis matrix of full column-rank. (Use of the symbol $\nu$ is to reflect the degrees of freedom in a standard test of such a hypothesis in the analogous univariate situation.)

After some initial preprocessing, functional data may typically be represented by a discrete, high-dimensional model. Such representations are made here using Fourier decomposition, whose advantage in FDA is not only to discretize the model but also to decorrelate the error structure and offer meaningful descriptive summarization. This is demonstrated on an existing data set in Section 2. (See also Fan and Lin, 1998, Spitzner, Marron, and Essick, 1998, Spitzner and Woodall, 2003, and Spitzner, 2008B for further demonstrations in FDA.) Also laid out in that section is how, taking into consideration the specific $\boldsymbol{L}$, data collected under the model (1.1) may be translated to that of the model

$$\boldsymbol{Y}_{n,j} = \boldsymbol{\theta}_j + n^{-1/2}\mathbf{e}_{n,j}, \quad (1.2)$$

for $j = 1, \ldots, p_n$, where $p_n$ represents some (high) maximum number of dimensions to be accounted for at a given $n$. The statistics $\boldsymbol{Y}_{n,j} = [Y_{n,j1}, \ldots, Y_{n,j\nu}]^T$ are $\nu$-dimensional discrete data vectors, each entry of which a linear combination Fourier coefficients (defined by a distinct column of $\boldsymbol{L}$) that partly summarizes the information relevant to the specific null hypothesis. The $\boldsymbol{\theta}_j = [\theta_{j1}, \ldots, \theta_{j\nu}]^T$ are mean vectors, and the $\mathbf{e}_{n,j} = [e_{n,j1}, \ldots, e_{n,j\nu}]^T$ are zero-mean, unit-covariance error-vectors such that the $e_{n,jk}$ are independent across $k$. (Across $j$, however, small correlations among the $e_{n,jk}$ are possible.) The functional linear hypothesis translates to

$$\mathrm{H}_0 : \boldsymbol{\theta}_j = 0 \text{ for } j = 1, \ldots, p_n \text{ \textit{versus} } \mathrm{H}_1 : \text{ not } \mathrm{H}_0, \tag{1.3}$$

in which the matrix $\boldsymbol{L}$ has been absorbed into the transformation from (1.1) to (1.2).

A typical assumption made in FDA is that the functional parameter $\boldsymbol{\beta}(t)$ in (1.1) is somehow "smooth." From an intuitive standpoint, this means that the $\beta_1(t), \ldots, \beta_P(t)$ are each taken to be a conglomeration of mainly large-scale, sweeping shapes, which are represented by low-frequency Fourier components. For rigorous analysis, smoothness is expressed more technically in Section 3 by restricting $\boldsymbol{\beta}(t)$ to a Sobolev class. At any rate, a key issue in testing is how to exploit the smoothness assumption, so as not to waste statistical power attempting to distinguish the "rougher" aspects of the model (*i.e.*, the small-scale wiggly shapes). This is especially important for testing in high-dimensions, where the vastness of the parameter space requires a careful management of power.

The transformation to the discrete model (1.2) typically assigns smaller $j$ to model-components associated with smoother functional attributes. Noting this, one would want a test that focuses power primarily on the model's lower-indexed components. One class of tests that do this is defined by test statistics of the form

$$Q_n = n \sum_{j=1}^{p_n} w_{n,j} \|\boldsymbol{Y}_{n,j}\|^2, \tag{1.4}$$

where each $0 < w_{n,j} \leq 1$ and $w_{n,j} \to 0$ as $j \to \infty$ to emphasize the $\boldsymbol{Y}_{n,j}$ with smaller $j$. The particular test of interest in this article is defined by the weight setting $w_{n,j} = j^{-1/2}$, which rewrites (1.4) as $Q_n^{OPT} = n \sum_{j=1}^{p_n} j^{-1/2} \|Y_{n,j}\|^2$. This technique of managing power by direct down-weighting, as in (1.4), shall be referred to as "tapering." Detailed asymptotic power properties of tests based on tapering are deduced in Spitzner (2008A), some of whose results are reproduced in Section 3. There, it is shown that $Q_n^{OPT}$ manages asymptotic power in an optimal way among test statistics of the form (1.4).

The asymptotic performance criteria used in Spitzner (2008A) and adopted here are taken from "rates of testing" theory, a framework articulating the rate at which power is retained in high-dimensional testing problems under geometric smoothness constraints. Its basic components are laid out in Ingster (1993), Spokoiny (1996), Lepski and Spokoiny (1999), Horowitz and Spokoiny (2001), and Gayraud and Pouet (2005), among others. Relevant criteria are described fully in Section 3. Within this context, the problem of selecting the $w_{n,j}$ in (1.4) for good asymptotic power would be aptly described as constrained rate-optimization among tests based on tapering, for which $Q_n^{OPT}$ is a solution.

There are, however, two types of asymptotic optimality to consider, and it is known that any test based on tapering is suboptimal with respect to the type that is more relevant for consideration in practice. Specifically, Ingster (1993) and Spokoiny (1996) respectively deduce "minimax" and "adaptive-minimax" rates, which bound the performance of any test with respect to an adopted smoothness geometry. The former rate defines a stronger form of asymptotic optimality, but the latter "adaptive" type is the one more practically relevant. It is well established that tests based on tapering can achieve Ingster's minimax rate of testing (*e.g.,* Ingster (1993) and Fan, Zhang, and Zhang (2001) provide two distinct examples), but Spitzner (2008A) deduces the suboptimality of any test based on tapering, including that based on $Q_n^{OPT}$, with respect to Spokoiny's adaptive-minimax criterion.

Existing tests that are known to achieve adaptive minimaxity are based on the alternative test-construction techniques of "truncation" or "thresholding." Spokoiny (1996) provides an adaptive-minimax test based on thresholding, and Fan, Zhang, and Zhang (2001) establish the adaptive minimaxity of Fan's (1996) "adaptive Neyman test," a test based on truncation. Details of these tests are provided in Section 4. Both have been developed for use in FDA in Fan and Lin (1998) and Abramovich *et al.* (2002).

Despite its asymptotic suboptimality, the test based on $Q_n^{OPT}$ is demonstrated in Section 4 to retain good and even superior power over the adaptive-minimax tests above in non-asymptotic, practically realistic FDA settings, illustrating that the improvements offered by thresholding or truncation, though guaranteed, may arise quite slowly asymptotically. That is, the evidence of this article suggests that, among tapering, truncation, and thresholding, the most powerful tests in typical FDA applications are constructed by tapering!

A novel aspect of this article is that it highlights a distinction between the asymptotic setup of FDA and that of other high-dimensional problems involving smoothness constraints, such as goodness-of-fit testing, within which many of the existing high-dimensional tests were first developed. The particular distinction has to do with the dependence of the dimensionality parameter, $p_n$, on the sample-size parameter $n$. In non-FDA scenarios, there is typically an explicit connection between these two parameters (often $p_n = n$), whereas in FDA the connection is largely hypothetical, and one typically has little or no control over the rate at which $p_n$ increases. Accordingly, an important concern in FDA is the sensitivity of test performance on the rate of $p_n \to \infty$. It shall be seen that the optimality property of the test based on $Q_n^{OPT}$ is robust in this regard.

### *1.1. Organization*

This paper is organized as follows. Section 2 presents an applied data example, which demonstrates the transformation from (1.1) to (1.2) and the test based on $Q_n^{OPT}$. This will provide grounding and intuition for subsequent discussion. The setup and relevance of asymptotic analysis in FDA is discussed in that section as well. Section 3 defines rates-of-testing criteria and presents the paper's main theoretical results. Section 4 gives details of several existing high-dimensional tests, and reports on empirical comparisons carried out by simulation. Concluding discussion appears in Section 5, in which the importance of studying tests based on tapering in FDA is further elaborated.

## 2. Functional data analysis and its asymptotic framework

Let us begin discussion with an analysis of the Canadian temperature data of Ramsay and Silverman (2005, ch. 13) under a functional linear model. (The data are available in a supplemental website to the book.) Subsequently, a general asymptotic setup for functional data analysis will be laid out.

### *2.1. An example data set*

The Canadian temperature data consist of daily mean-temperature profiles across the year at 31 weather stations in three regions of Canada: there are $M_1 = 14$ "Atlantic" stations, $M_2 = 5$ "Pacific" stations and $M_3 = 12$ "Continental" stations. The raw measurements are displayed in the top panels of Figure 1. (The original data set analyzed by Ramsay and Silverman also includes three stations in an additional "Arctic" region, and one additional Atlantic station, at "Schefferville." For the present analysis it makes sense to have set aside the Arctic stations since there are so few of them, and the Schefferville station since its location is of unusually high latitude relative to the other Atlantic stations.) Ramsay and Silverman remark that the region-effects seen in these data are "more complex than the constant or even sinusoidal effects that one might expect," but note specifically that the Pacific stations tend to have warmer winter temperatures and Continental stations tend to have colder winter temperatures, while all regions' summer temperatures tend close to the average. The latter observations refer to large-scale attributes of the temperature profiles, and it will be presumed that these and other large-scale attributes are of primary interest. Nevertheless, smaller-scale, "more complex" attributes are not to be ignored; we want a systematic way to explore essentially all aspects of the data. Fourier decomposition provides just such a technique, one that is furthermore most appropriate for data such as yearly temperature profiles whose attributes tend to be periodic.

### *2.2. Fourier decomposition*

The notation used in this analysis will parallel that of Section 1, but, for simplicity, without the subscript $n$, thereby ignoring any replication concepts until further ideas are laid down. Set $M = M_1 + M_2 + M_3 = 31$ and denote by $Y_k(t)$ the measurement of the $k$'th station at time $t$ within a typical year: the measurement times are $t_1, \ldots, t_r$, where $t_1$ is January 1, $r = 365$, and $t_l = l$, so that $t_{l+1} - t_l$ is one day. Next define the Fourier basis functions on $(0, 365]$ according to $\psi_1(t) = 1$, $\psi_{2j}(t) = \sin(j\pi(2t/r - 1))$, and $\psi_{2j+1}(t) = \cos(j\pi(2t/r - 1))$ for $j = 1, 2, \ldots$ The $j$'th Fourier coefficient of the $k$'th station is then $Y_{jk}^* = r^{-1} \sum_{l=1}^{r} Y_k(t_l) \psi_j(t_l)$, which gives the $j$'th coefficient of a multiple regression of $Y_k(t_l)$ across $t_1, \ldots, t_r$ onto any finite set of regressors $\psi_j(t_l)$ across $t_1, \ldots, t_r$. Coefficients associated with the first few $\psi_j$ are shown in the bottom panels of Figure 1, centered and scaled for each $j$ so that the mean and sum of squares of the displayed $Y_{jk}^*$ match common values.

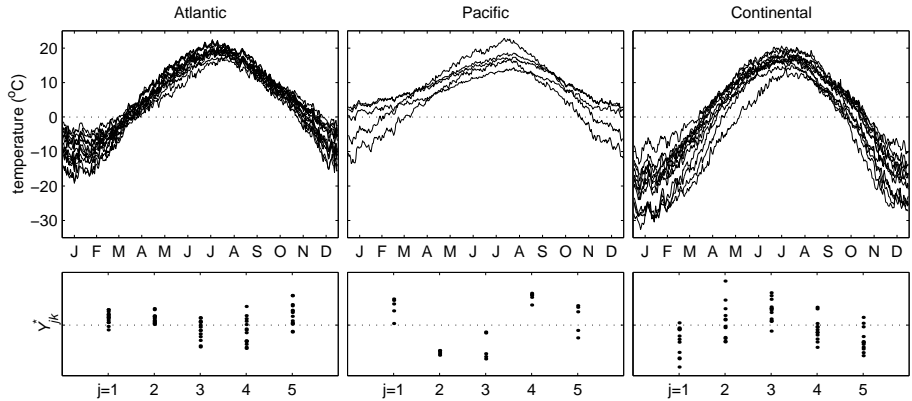

FIG 1. *Canadian temperature data for weather stations in the Atlantic, Pacific, and Continental regions. The top panels plot average temperatures $Y_k(t)$ across the year, and the bottom panels plot individual sets of Fourier coefficients $Y_{jk}^*$ for $j = 1, \ldots, 5$.*

Consider, for instance, that the shapes of $\psi_1(t) = 1$ and $\psi_3(t) = \cos(\pi(2t/r - 1))$ convey interpretations whereby $Y_{1k}^*$ measures the $k$'th station's yearly average temperature and $Y_{3k}^*$ measures its differential between the winter and summer temperatures. With this in mind, observe from the bottom panels of Figure 1 that the $Y_{1k}^*$ tend to be larger for the Pacific region and smaller for the Continental region, while the $Y_{3k}^*$ tend to be smaller for the Pacific region and larger for the Continental region. This reflects the observations made by Ramsay and Silverman: between these two regions, the yearly average temperature of the Pacific region is warmer (as reflected in the $Y_{1k}^*$) and there is a smaller differential between the winter and summer temperatures (as reflected in the $Y_{3k}^*$). Each of the remaining sets of coefficients describe a distinct attribute of these data. For instance, the $Y_{2k}^*$ describe asymmetries between the spring and fall transition periods, and the $Y_{jk}^*$ with larger $j$ summarize finer periodic attributes of the stations' yearly profiles.

Fourier decomposition of these data is not discussed in Ramsay and Silverman (2005) itself, but it is in the book's supplemental internet materials, within which there appears the remark: "it was decided that 65 basis functions captured enough of the detail in the temperature data...," referring to Fourier basis functions. The present analysis will follow this guideline and consider just the Fourier coefficients $Y_{jk}^*$ with $j = 1, \ldots, 65$.

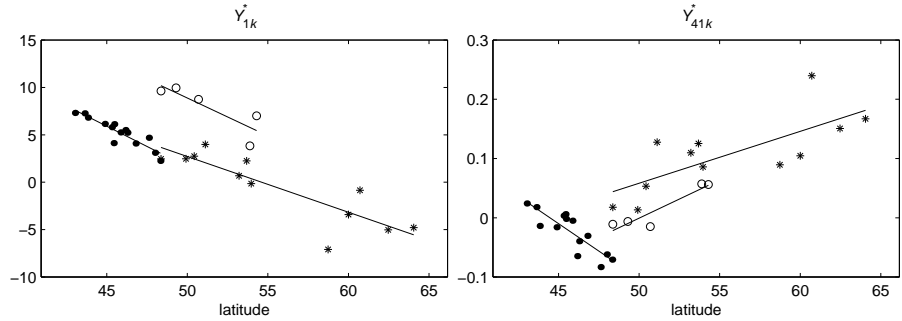

FIG 2. *Scatterplots of $Y^*_{1k}$ and $Y^*_{41k}$ by latitude, with region-specific fitted regression lines. The Atlantic stations are indicated by dots (●), the Pacific stations by circles (○), and the Continental stations by asterisks (∗).*

### 2.3. Linear models and testing

Ramsay and Silverman (2005) treat the temperature profiles by an elementary multi-group functional model, $Y_{k(g,l)}(t) = \mu_g(t) + \epsilon_{k(g,j)}(t)$, where $k = k(g,l)$ indexes the $l$'th station of region $g$, $\mu_g(t)$ is the mean temperature profile for the $g$'th region, and $\epsilon_{k(g,j)}(t)$ are random errors. Their analysis concludes that there are indeed vast differences in temperature profiles among the regions. Consider, however, that one would expect latitude to explain some substantial portion of the variation in temperature from station to station, and moreover that it is surely possible for the dependency of temperature on latitude to change from region to region. Accordingly, the model considered here is an extension of the multi-group model that incorporates latitude as a covariate in such a way as to allow for differences in both intercept and slope across regions. Denoting by $x_k$ the latitude of the $k$'th station, and writing $\bar{x} = M^{-1}\sum_{k=1}^{M} x_k$, the extended model is $Y_{k(g,l)}(t) = \mu_g(t) + \beta_g(t)(x_{k(g,j)} - \bar{x}) + \epsilon_{k(g,j)}(t)$, a functional linear model with $P = 6$ regression parameters, $\mu_g(t)$ and $\beta_g(t)$ for $g = 1, 2, 3$. From this, an analogous model is implied for each set of Fourier coefficients, which for the $j$'th set is

$$Y^*_{k(g,l)j} = \mu_{gj} + \beta_{gj}(x_{k(g,l)} - \bar{x}) + \epsilon_{k(g,l)j}. \tag{2.1}$$

Component-specific estimates and tests may be carried out under the model (2.1) using standard linear-regression methodology. To illustrate, Figure 2 displays scatterplots of Fourier coefficients corresponding to $j = 1$ in the left panel and $j = 41$ in the right panel, with estimated region-specific fitted regression lines drawn in. The former panel depicts patterns one would more-or-less expect to see in these data. There, the $Y^*_{1k}$ are seen to generally decrease as latitude increases, as would be expected of measurements of yearly average-temperature. Observe also that the $Y^*_{1k}$ of the Atlantic and Continental stations follow a common trend fairly consistently, whereas those of the Pacific stations fall above that trend. The fitted slopes of all three regions are roughly the same. A very different pattern is seen in the scatterplot of the $Y^*_{41k}$. Among those coefficients it is seen that the Pacific and Continental stations follow roughly the same trends, though shifted slightly, and that those trends are very different from the trend followed by the Atlantic stations: the $Y^*_{41k}$ tend to increase with latitude for the former stations but decrease for the latter.

TABLE 1
*P-values for component-specific and global tests comparing the relationships between temperature and latitude across regions. The headings "falls on common trend" and "same slope" refers to the respective null hypotheses $H_0 : (\mu_{g_1 j} - \bar{\mu}_j) + (\beta_{g_1 j} - \bar{\beta}_j)(\bar{x}_{g_1} - \bar{x}) = (\mu_{g_2 j} - \bar{\mu}_j) + (\beta_{g_2 j} - \bar{\beta}_j)(\bar{x}_{g_2} - \bar{x})$ and $H_0 : \beta_{g_1 j} = \beta_{g_2 j}$, where $(g_1, g_2)$ is taken across all pairs of distinct regions consistent with the row label. The headings $Y^*_{1k}$ and $Y^*_{41k}$ indicate component-specific tests at $j = 1$ and $j = 41$, respectively, each with $\nu$ numerator degrees of freedom. "Global" p-values are simulated from the null distribution of $F_{global}$.*

| Comparison | $\nu$ | Falls on common trend | | | Same slope | | |
|---|---|---|---|---|---|---|---|
| | | $Y^*_{1k}$ | $Y^*_{41k}$ | Global | $Y^*_{1k}$ | $Y^*_{41k}$ | Global |
| All regions | 2 | <0.001 | <0.001 | <0.001 | 0.480 | <0.001 | 0.005 |
| Atlantic to Pacific | 1 | <0.001 | 0.019 | <0.001 | 0.881 | <0.001 | 0.031 |
| Atlantic to Continental | 1 | 0.289 | <0.001 | <0.001 | 0.313 | <0.001 | 0.003 |
| Pacific to Continental | 1 | <0.001 | <0.001 | <0.001 | 0.454 | 0.510 | 0.180 |

The pattern associated with the $Y^*_{41k}$ within a yearly temperature profile is symmetric with a period about two-and-a-half weeks. Its magnitude is quite small, accounting for relative temperature differences of less than one half of one degree. Yet our analysis suggests its relationship with latitude may distinguish the Atlantic region from the others.

Formal tests of the patterns noted above may be formulated in terms of the parameters of the model (2.1). Corresponding p-values are listed in Table 1. To describe the relevant hypotheses, it will be convenient to define an "overall" regression line at each $j$, $\overline{REG}_j(x) = \bar{\mu}_j + \bar{\beta}_j(x - \bar{x})$, where $\bar{\mu}_j = M^{-1} \sum_g M_g^3{}_{g=1} \mu_{gj}$ and $\bar{\beta}_j = M^{-1} \sum_{g=1}^3 M_g \beta_{gj}$. Referring to the headings of Table 1, the null hypothesis labeled "falls on common trend" is $H_0 : (\mu_{g_1 j} - \bar{\mu}_j) + (\beta_{g_1 j} - \bar{\beta}_j)(\bar{x}_{g_1} - \bar{x}) = (\mu_{g_2 j} - \bar{\mu}_j) + (\beta_{g_2 j} - \bar{\beta}_j)(\bar{x}_{g_2} - \bar{x})$ for $(g_1, g_2) \in G$, for which $\bar{x}_g = M_g^{-1} \sum_{l=1}^{M_g} x_{k(g,l)}$, where the index-set $G$ consists of all pairs of distinct indices among regions indicated in the table's corresponding row-label. The null hypothesis labeled "same slope" is $H_0 : \beta_{g_1 j} = \beta_{g_2 j}$ for $(g_1, g_2) \in G$.

Each null hypothesis above may be written as a linear hypothesis $H_0 : \boldsymbol{L}^T \boldsymbol{\beta}_j = 0$, where $\boldsymbol{\beta}_j = [\mu_{1j}, \mu_{2j}, \mu_{3j}, \beta_{1j}, \beta_{2j}, \beta_{3j}]^T$ and $\boldsymbol{L}$ is a matrix of full column-rank that is determined by the specific hypothesis. The associated test statistic is a ratio of independent "mean-square" statistics $F_j = \text{MS}_j(\boldsymbol{L})/\text{MSE}_j$, which follows an F distribution whose non-centrality parameter is zero only under $H_0$, assuming the errors in (2.1) are independent, homoscedastic, and Gaussian. (Explicit formulas are provided in Section 2.4; see also, Seber, 2003.) Associated degrees of freedom are $\nu = \text{rank } \boldsymbol{L}$, whose values are indicated in Table 1, and $M - P = 25$.

The p-values in Table 1 reflect the observations made above on Figure 2. Comparing against the standard 0.05 level, the p-values in the column labeled $Y^*_{1k}$ under "falls on common trend" separate out the Pacific region as falling off of a common trend that is followed by the other regions. Those in the $Y^*_{1k}$ column under "same slope" indicate no evidence for differences in slope. Similarly, the column labeled $Y^*_{41k}$ under "same slope" separates out the Atlantic region as having a different relationship with latitude than the others.

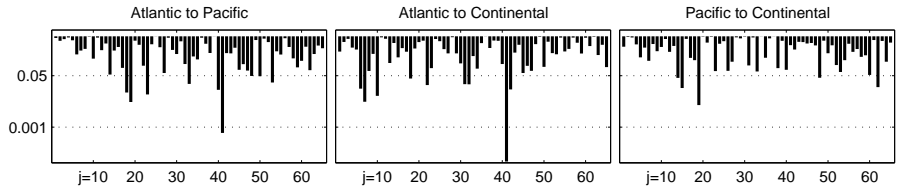

Fig 3. *P-values for component-specific tests of "same slope." The null hypotheses are $H_0 : \beta_{g_1,j} = \beta_{g_2,j}$, for region indices $(g_1, g_2)$ indicated above each plot. The vertical axis is on a logarithmic scale; reference levels are indicated at $\alpha = 0.001$ and $\alpha = 0.05$.*

There remains the question of how to carry out these tests "globally." That is, can the same hypotheses be tested on the whole of the temperature profiles, inasmuch as they are represented by 65 sets of Fourier coefficients? This type of question is the central concern of this article. A starting point is to consider the plots in Figure 3, which charts the p-values from component-specific tests of "same slope" across all sets of Fourier coefficients, $j = 1, \ldots, 65$, each panel corresponding to a separate pair of regions. The most significant differences in slope are reflected in the very small p-values displayed in the two leftmost panels, at $j = 41$. Other p-values are "small" as well, in the sense of falling below the 0.05 level, even in the rightmost panel, but nowhere near as small as these two at $j = 41$. However, in light of there being 65 test results to examine per panel, it is no surprise to find at least a handful of small p-values. Of interest, then, is to deduce a single assessment for each panel which combines the p-values across all $j = 1, \ldots, 65$.

Well suited to this task is a test statistic constructed by tapering, in a manner similar to (1.4). For the present situation, let us define this statistic as $F_{global} = \sum_{j=1}^{65} w_j F_j$, and set the weights to $w_j = j^{-1/2}$, paralleling the construction of $Q_n^{OPT}$. By combining the individual test statistics $F_j$ this way, $F_{global}$ is globally sensitive, but it down-weights the influence of the $Y^*_{jk}$ with larger $j$, as is desired to reflect primary interest in the larger-scale shapes.

Independence shall be assumed among the $Y^*_{jk}$ across $j$ (which is justified in in Section 2.4), so that the null distribution of $F_{global}$ is fully defined. Simulated p-values (using one million iterations) for global versions of the null hypotheses discussed above are listed in the columns of Table 1 labeled "global." Comparing against 0.05, those in the "same slope" portion of the table again separate out the Atlantic region as having a different relationship with latitude than the others, but this conclusion now accounts for essentially all aspects of the temperature profiles.

### 2.4. An asymptotic framework for functional data analysis

The sort of analysis carried out on the Canadian temperature data is now cast in the general context described in Section 1. The transformation from the continuous model (1.1) to the discrete model (1.2) will be laid out in detail, and its properties discussed. Hereafter, the notation will revert back to that of the replicated model, in which there is an explicit reference to a sample-size parameter, $n$.

Regarding the setup for asymptotic analysis, there is some imprecision in defining replication strictly according to the model (1.1), for the presence of a continuous covariate would make it unrealistic to expect exact duplication of the regressor matrix, $\boldsymbol{X}$, except in carefully designed experiments. Absent a continuous covariate, however, many multi-group experiments do allow the parameterization $M_g = N_g n$, where $M_g$ are group-specific sample sizes, and the "sample size" $n$ is some common divisor among them. In such cases the matrix $\boldsymbol{X}$ would consist entirely of zeros and ones, and $N = \sum_g N_g$. The notation adopted in the Canadian temperature example is intended to suggest such a multi-group experiment; yet, the corresponding replication concept is unrealistic, for if any additional weather stations were sampled, one would not expect the latitudes of those stations to match any of those in the current data set.

Such complications notwithstanding, the purpose of introducing the sample-size parameter, $n$, is to manifest the notion that an increase in the amount of data collected is coupled with a decrease in error variability. This is apparent in expression (1.2), which shows the magnitude of error in the discrete model to shrink at the rate $n^{-1/2}$. Some readers might prefer to reinterpret the asymptotic formulation used in this article to one in terms of shrinking errors, in which case the results presented here would directly translate. Otherwise, a "pure" interpretation of replication in the model (1.1) might take $n = 1$ in an analysis of current data, and treat any future replication as entirely hypothetical; or, one might employ various conceptual devices, such as sampling covariates from a distribution, to modify the model (1.1) and its translation to (1.2) so as to make replication more realistic. Despite these possibilities, the perspective taken here is that the model (1.1) is entirely adequate for illustrating the key ideas of present interest, and any modification would only add technical complications that are tangential to them.

Regarding data collection and Fourier decomposition, the functional measurements are assumed to have been taken along a dense, finite grid that is common to all $\boldsymbol{Y}_i(t)$, as in the Canadian temperature data. For $t$ in the domain $(a, b]$, the points of the grid are taken to be $t_l = a + (b-a)l/r$ for $l = 1, \ldots, r$ and some fixed, large $r$ (which may not be $p_n$). The data associated with the curve $Y_{i,k}(t)$ are $Y_{i,k}(t_1), \ldots, Y_{i,k}(t_r)$. (In more general situations the grid may change from curve to curve, but to avoid additional complication it will be assumed a good approximation to the present setup is available, *e.g.*, by interpolating measurements onto a fixed grid.) Set $\psi_1(t) = 1$, $\psi_{2j}(t) = \sin(\pi j\{2(t-a)/(b-a) - 1\})$, and $\psi_{2j+1}(t) = \cos(\pi j\{2(t-a)/(b-a) - 1\})$ for $j = 1, 2, \ldots$ Writing $\boldsymbol{Y}^*_{i,j} = [Y^*_{i,j1}, \ldots, Y^*_{i,jN}]^T$, the Fourier coefficients, $Y^*_{i,jk}$, are calculated by the formula

$$\boldsymbol{Y}^*_{i,j} = \frac{1}{r}\sum_{l=1}^{r}\boldsymbol{Y}_i(t_l)\psi_j(t_l). \quad (2.2)$$

To make the final step to the discrete model (1.2), a linear transformation involving $\boldsymbol{L}$ is used to tailor the statistics (2.2) to the linear hypothesis. Set $\boldsymbol{H} = \{\boldsymbol{L}^T(\boldsymbol{X}^T\boldsymbol{X})^{-1}\boldsymbol{L}\}^{-1/2}$ and define the $\boldsymbol{Y}_{n,j}$ in (1.2) according to

$$\boldsymbol{Y}_{n,j} = \frac{1}{\sigma_j n}\sum_{i=1}^{n}\boldsymbol{H}\boldsymbol{L}^T(\boldsymbol{X}^T\boldsymbol{X})^{-1}\boldsymbol{X}^T\boldsymbol{Y}^*_{i,j}, \quad (2.3)$$

where $\sigma_j^2 = V[Y^*_{i,jk}]$. The remaining objects defining (1.2) are

$$\boldsymbol{\theta}_j = \frac{1}{\sigma_j r}\sum_{l=1}^{r}\boldsymbol{H}\boldsymbol{L}^T\boldsymbol{\beta}(t_l)\psi_j(t_l) \text{ and} \quad (2.4)$$

$$\boldsymbol{e}_{n,j} = \frac{1}{\sigma_j r\sqrt{n}}\sum_{i=1}^{n}\sum_{l=1}^{r}\boldsymbol{H}\boldsymbol{L}^T(\boldsymbol{X}^T\boldsymbol{X})^{-1}\boldsymbol{X}^T\boldsymbol{\epsilon}_i(t_l)\psi_j(t_l).$$

An often-appropriate assumption has each $\boldsymbol{\epsilon}_i$ a stationary process such that $\epsilon_{i,k}(t_l) = \sum_{m=-\infty}^{\infty}\gamma_m\ \eta_{i,k}(t_l - (b-a)m/r)$, for which $\eta_{i,k}(a + (b-a)m/r)$ is, across integer $m$, a mean-zero independent and identically distributed sequence with finite fourth moment, and $\sum_{m=-\infty}^{\infty}|\gamma_m| < \infty$. When this assumption is valid, Theorem 10.3.2.*i* of Brockwell and Davis (1991) implies that each $Cov(Y^{*2}_{i,jl}, Y^{*2}_{i,kl}) \to 0$ as $r \to \infty$, for $j \neq k$. (See also Corollary 3.1.1.*i* in Section 3, below.) Thus, Fourier decomposition provides a means to decorrelate the functional linear model, while the statistics (2.2) or (2.3) capture its core structure.

In the typical case where the $\sigma_j^2$ are unknown, the test statistic $Q_n$ in (1.4) would be replaced by $\hat{Q}_n = n\sum_{j=1}^{p_n}w_j\|\hat{\boldsymbol{Y}}_{n,j}\|^2$, where $\hat{\boldsymbol{Y}}_{n,j}$ is defined as in the right side of (2.3) but with an estimate $\hat{\sigma}^2_{n,j}$ substituting for $\sigma_j^2$. For instance, the $\hat{\sigma}^2_{n,j}$ may be the usual unbiased estimates

$$\hat{\sigma}^2_{n,j} = \frac{1}{(nN-P)}\sum_{i=1}^{n}\sum_{k=1}^{N}(Y^*_{i,jk} - \boldsymbol{x}_k^T\hat{\boldsymbol{\beta}})^2, \quad (2.5)$$

where

$$\hat{\boldsymbol{\beta}} = \frac{1}{n}\sum_{i=1}^{n}(\boldsymbol{X}^T\boldsymbol{X})^{-1}\boldsymbol{X}^T\boldsymbol{Y}^*_{i,j}.$$

This modification led to the statistic $F_{global}$ in the Canadian temperature example, for which $F_j = \text{MS}_j(\boldsymbol{L})/\text{MSE}_j$ has

$$\text{MS}_j(\boldsymbol{L}) = n^{-1}\sum_{i=1}^{n}\boldsymbol{H}\boldsymbol{L}^T(\boldsymbol{X}^T\boldsymbol{X})^{-1}\boldsymbol{X}^T\boldsymbol{Y}^*_{i,j}/\nu \quad \text{and} \quad \text{MSE}_j = \hat{\sigma}^2_{n,j},$$

so that $F_{global} = \sum_{j=1}^{p_n} w_j F_j = \hat{Q}_n/\nu$. Fan and Lin (1998, sec. 3.4) discuss other estimates of $\sigma_j^2$ based on smoothing across $j$, which are specialized for the Fourier decomposition technique.

Regarding the relationship between $p_n$ and $n$, the reader should notice that $p_n$ never appears in the transformation from (1.1) to (1.2), and so these two parameters are never actually connected. The parameter $p_n$ is constrained by the resolution of the grid $t_1, \ldots, t_r$, which to avoid numerical error requires $p_n \leq r$. Moreover, it is possible that $p_n$ would be set according to subjective modeling assumptions such as an analyst's determination of the number of basis functions needed to describe the data, as was made in the Canadian temperature example. In Spitzner, Marron, and Essick (1998), $p_n$ is set subjectively to avoid observed defects in the ability of the $\psi_j$ to decorrelate the model at larger $j$. At any rate, the justification for taking $n \to \infty$ and $p_n \to \infty$ is that it forms an appropriate abstract conceptualization for repeated measurement of functional data in accordance with a *global* point of view. In particular, $p_n \to \infty$ represents a situation where the grid $t_1, \ldots, t_r$ is to become increasingly dense, and if $n \to \infty$ as well then potentially all available information about the curve model will be captured in the limit.

## 3. Rates of testing for the tapering mechanism

The discussion now turns to rates-of-testing theory and the asymptotic optimality of $Q_n^{OPT}$ among tests based on tapering. The more technical aspects of this discussion have been omitted, but can be found in Spitzner (2008A).

Smoothness constraints are formally defined within rates-of-testing theory as a restriction of the functional parameter $\boldsymbol{\beta}(t)$ in (1.1) to a smooth-function class. In the most general settings, this would be a Besov class, but here it is taken to be a Sobolev class, a special case, which is appropriate when working with Fourier decompositions. Such constraints may be expressed as a restriction of the mean vectors of the discrete model (1.2) to the geometry

$$\mathcal{B}_{s,M} = \left\{ (\boldsymbol{\theta}_1, \boldsymbol{\theta}_2, \ldots) : \sqrt{\sum_{j=1}^{\infty} j^{2s} \|\boldsymbol{\theta}_j\|^2} \leq M \right\}, \quad (3.1)$$

a Sobolev ellipsoid of radius $M$ in infinite-dimensional discrete space, where $M > 0$ and $s > 1/2$ are fixed constants. The notation $\tilde{s} = 4s + 1 > 3$ shall also be used. The bound on the norm in (3.1) models smoothness by restricting expression of the higher-indexed $\boldsymbol{\theta}_j$, with larger $s$ making the restriction stronger. Moreover, Parseval's identity implies that $(\boldsymbol{\theta}_1, \boldsymbol{\theta}_2, \ldots) \in \mathcal{B}_{s,M}$ is equivalent to the assumption that the corresponding $\boldsymbol{\beta}(t)$ in (1.1), assuming the $\boldsymbol{\theta}_j$ arise through (2.4), is an element of a Sobolev ellipsoid in continuous space, $\{\boldsymbol{\beta}(t) = [\beta_1(t), \ldots, \beta_P(t)]^T : \| \int_0^1 \beta_k^{(s_c)}(t)dt \| \leq M_c\}$ for some $s_c = 1, 2, \ldots$ and $M_c > 0$, which are easily determined. (For details and further discussion, see Adams and Fournier, 2003).

The "rates" in rates-of-testing theory, which characterize test performance, are described as follows. Fix $s > 1/2$, $M > 0$, and for each $n$ let $\phi_n = \phi_n(\boldsymbol{Y}_{n,1}, \ldots, \boldsymbol{Y}_{n,p_n})$ be a test of (1.3) such that $\lim_n \mathrm{P}_0[\phi_n = 1] \leq \alpha$, for a fixed level $\alpha \in (0,1)$. The notation $\mathrm{P}_{\boldsymbol{\theta}}$ is here used to denote probabilities under the model (1.2) for a specific $\boldsymbol{\theta} = (\boldsymbol{\theta}_1, \boldsymbol{\theta}_2, \ldots)$ and fixed $n$. Rates-of-testing criteria are formulated from sequences $\delta_n \to 0$ satisfying

$$\inf_{\boldsymbol{\theta} \in \mathrm{H}_1(\delta_n/\delta_n^*; s, M)} \mathrm{P}_{\boldsymbol{\theta}}[\phi_n = 1] \to 1 \text{ for every } \delta_n^* \to 0, \tag{3.2}$$

where

$$\mathrm{H}_1(\delta; s, M) = \left\{ \boldsymbol{\theta} \in \mathcal{B}_{s,M} : \sqrt{\sum_{j=1}^{\infty} \|\boldsymbol{\theta}_j\|^2} \geq \delta \right\}. \tag{3.3}$$

The criterion (3.2) describes the rate at which a gap may shrink between the null hypothesis and a class of "distinguishable" alternatives, those the test would be able to detect with high power, asymptotically. The better-performing tests allow this gap to shrink faster: if for some $\delta_n \to 0$ the criterion (3.2) is satisfied for one test, but not another, the former test is preferred.

Ingster's (1993) minimax performance bound states that for no test does any $\delta_n = o(n^{-2s/\tilde{s}})$ satisfy (3.2), but there is a test (based on tapering) for which $\delta_n = n^{-2s/\tilde{s}}$ satisfies (3.2). This identifies the rate $\hat{\delta}_n^M(s) = n^{-2s/\tilde{s}}$ as minimax for the geometry $\mathcal{B}_{s,M}$ at a specific $s$. Suppose now that fixed bounds $s_* < s^*$ are given, and for each $s_* < s < s^*$ one is to consider a separate sequence $(\delta_n(s))$, and set $\hat{\delta}_n^{AM}(s) = \{n^2(\log \log n)^{-1}\}^{-s/\tilde{s}}$. Spokoiny (1996) establishes that for no test is (3.2) satisfied across $s_* < s < s^*$ if $\delta_n(s) = o(\hat{\delta}_n^{AM}(s))$ for some such $s$. It is also shown there is a test (based on thresholding) for which $\delta_n(s) = \hat{\delta}_n^{AM}(s)$ does satisfy (3.2) across $s_* < s < s^*$. This identifies the rates $\hat{\delta}_n^{AM}(s)$ as adaptive-minimax for $\mathcal{B}_{s,M}$ across $s_* < s < s^*$. (The optimal tests alluded here are the same mentioned in Section 1 and are described later in Section 4.)

The main technical result for evaluating tests based on tapering rewrites the criterion (3.2) in terms of the parameters of the test statistic (1.4).

**Theorem 3.1.** *Assume the model (1.2) and suppose $(Q_n)$ is a sequence of test statistics with each $Q_n$ as in (1.4) for associated sequences $(w_{n,j})$ and $(p_n)$ such that each $0 < w_{n,j} \leq 1$ and $p_n \to \infty$ as $n \to \infty$. Set $S_n(p) = w_{n,1}^2 + \cdots + w_{n,p}^2$, $W_n(p) = \min\{w_{n,j}^2 : j \leq p\}$, $U_n(p,q) = qW_n(q)/S_n(p)$, and $U_n(p) = U_n(p,p)$. Suppose at each $n$ the $e_{n,jk}$ are independent across $k$, and*

$$\left| V\left[ n \sum_{j=1}^{p_n} w_j \|\boldsymbol{Y}_{n,j}\|^2 \right] - n^2 \sum_{j=1}^{p_n} w_{n,j}^2 V[\|\boldsymbol{Y}_{n,j}\|^2] \right| = O(S_n(p_n)). \tag{3.4}$$

*Suppose further that each $e_{n,jk}$ is such that $E[e_{n,jk}] = 0$, $V[e_{n,jk}] = 1$, $E[e_{n,jk}^4] \asymp 1$ uniformly across $j$, $k$, and $n$, and $P[e_{n,jk} \leq -t] > 0$ for each $t > 0$. Let* ($\delta_n$) *be some positive sequence for which* $\delta_n \to 0$. *Fix* $\alpha$ *and, for each* $n$, *let* $\phi_n^Q$ *denote the size-*$\alpha$ *test which rejects the null hypothesis in (1.3) when* $Q_n$ *exceeds some critical value. For fixed* $s > 1/2$, $M > 0$, *and* $\phi_n = \phi_n^Q$ *the criterion (3.2) holds if, and only if, both*

$$\text{(i.)} \quad \limsup_{n \to \infty} n^2 U_n(p_n) p_n^{-\tilde{s}} < \infty \quad \text{and} \tag{3.5}$$

$$\text{(ii.)} \quad \liminf_{n \to \infty} n^2 U_n(p_n, q_n) q_n^{-\tilde{s}} > 0,$$

*where* $q_n = \{\delta_n/M\}^{-1/s}$. *The same conclusion holds if the* $Y_{n,jk}$ *in* $Q_n$ *are replaced with* $Y_{n,jk}(1+o_p(1))$, *provided* $Cov(Y_{n,jk}^2(1+o_p(1)), Y_{n,jl}^2(1+o_p(1))) \to Cov(Y_{n,jk}^2, Y_{n,jl}^2)$ *for each* $j$, $k$, *and* $l$.

*Proof.* This is Theorem 1 of Spitzner (2008A). □

An important corollary of this theorem establishes its validity under assumptions on the covariance structure of the continuous model (1.1) that would typically be made in practice. The case of unknown variances is also treated.

**Corollary 3.1.1.** *Suppose the model (1.2) derives from the functional linear model (1.1) via the transformation (2.3), and assume the notation of Section 2.4. The conclusions of Theorem 3.1 hold true under the following statement* (i) *and remain true under statement* (ii).

(i) *Each* $\epsilon_i$ *in (1.1) is such that* $\epsilon_{i,k}(t_l) = \sum_{m=-\infty}^{\infty} \gamma_m \, \eta_{i,k}(t_l - (b-a)m/r)$, *for which the coefficients* $\gamma_m$ *satisfy* $\sum_{m=-\infty}^{\infty} |\gamma_m||m|^{1/2} < \infty$, *the* $\eta_{i,k}(a+(b-a)m/r)$ *are independent and identically distributed across integer* $m$ *with* $E[\{\eta_{i,k}(a+(b-a)m/r)\}^4] < \infty$, *and the weight sequence* $(w_{n,j})$ *does not increase in* $j$ *for each* $n$.

(ii) $Q_n$ *is replaced by* $\hat{Q}_n$, *substituting for* $\sigma_j^2$ *the usual unbiased estimates* $\hat{\sigma}_{n,j}^2$, *defined in (2.5).*

*Proof.* Under the conditions of statement $(i)$, and assuming $p_n \le r$, Theorem 10.3.2.*ii* of Brockwell and Davis (1991) implies that $p_n Cov(Y_{i,jl}^{*2}, Y_{i,kl}^{*2})$ is uniformly bounded across $j, k = 1, \ldots, p_n$ with $j \ne k$. From this comes the property $Cov(Y_{n,jl}^2, Y_{n,kl}^2) < B/(n^2 p_n)$ for some $B$ across $j, k = 1, \ldots, p_n$ with $j \ne k$. Thus, if $(w_{n,j})$ is as indicated, one has

$$n^2 \sum_{j=1}^{p_n-1} \sum_{k=j+1}^{p_n} \sum_{l=1}^{\nu} w_{n,j} w_{n,k} Cov(Y_{n,jl}^2, Y_{n,kl}^2)$$

$$\le \frac{B\nu}{p_n} \sum_{j=1}^{p_n-1} (p_n - j) w_{n,j}^2 \le B\nu S_n(p_n),$$

and so the criterion (3.4) is satisfied. Statement $(ii)$ is readily verified using the delta rule to show that $E[\hat{\sigma}_{n,j}^{-2}] \to \sigma_j^{-2}$, $Cov(\hat{\sigma}_{n,j}^{-2}, \hat{\sigma}_{n,k}^{-2}) \to 0$, and $Cov(Y_{n,jk}, \hat{\sigma}_{n,j}^{-2}) \to 0$, from which it follows that $Cov(\hat{Y}_{n,jk}^2, \hat{Y}_{n,jl}^2) \to Cov(Y_{n,jk}^2, Y_{n,jl}^2)$. □

In Spitzner (2008A), Theorem 3.1 is applied to characterize and bound the rates-of-testing performance of tests based on tapering. The key results most relevant to present purposes are summarized as follows.

**Corollary 3.1.2.** *Assume the notation and conditions of Theorem 3.1.*

(i) *Suppose (3.5.i) holds and $\hat{\delta}_n \to 0$ is such that $n^2 U_n(p_n, \hat{q}_n)\hat{q}_n^{-\tilde{s}} \asymp 1$, where $\hat{q}_n = \{\hat{\delta}_n/M\}^{-1/s}$. Then, the sequence $(\hat{\delta}_n)$ defines a "boundary rate" of the test $\phi_n^Q$ in the sense that for any sequence $(\delta_n)$ the criterion (3.2) holds if $\hat{\delta}_n = O(\delta_n)$ but not if $\delta_n = o(\hat{\delta}_n)$.*

(ii) *Set $\hat{\delta}_n^Q(s) = \{n^2(\log n)^{-1}\}^{-s/\tilde{s}}$ and suppose $1/2 < s_* < s^*$. For no test $\phi_n^Q$ is (3.2) satisfied for $\delta_n = o(\hat{\delta}_n^Q(s))$ across $s_* < s < s^*$. Moreover, if $w_{n,j} = j^{-1/2}$, as in $Q_n^{OPT}$, and $p_n$ is such that $\{n^2(\log n)^{-1}\}^{1/3}/p_n$ is bounded and $\log p_n \asymp \log n$, then $\phi_n^Q$ satisfies (3.2) for $\delta_n = \hat{\delta}_n^Q(s)$ across $s_* < s < s^*$. Thus, $\hat{\delta}_n^Q(s)$ is an "optimal adaptive rate of testing for the tapering mechanism."*

*Proof.* Statement (*i*) and the first assertion of (*ii*) follow immediately from Theorems 2 and 3, respectively, of Spitzner (2008A). To prove the second assertion of (*ii*), first observe that the specified settings imply $\{n^2(\log n)^{-1}\}^{1/\tilde{s}}/p_n$ is bounded, since $\tilde{s} > 3$, and $S_n(p_n) \asymp \log p_n$. It follows that the sequence in (3.5.i) has

$$n^2 U_n(p_n)p_n^{-\tilde{s}} \asymp n^2 p_n^{-\tilde{s}}/\log p_n = \{(\log n)/(\log p_n)\} \left[\{n^2(\log n)^{-1}\}^{1/\tilde{s}}/p_n\right]^{\tilde{s}},$$

which is bounded. Next set $\delta_n = \hat{\delta}_n^Q(s)$, so that $q_n = M^{1/s}\{n^2(\log n)^{-1}\}^{1/\tilde{s}}$, and observe

$$n^2 U_n(p_n, q_n)q_n^{-\tilde{s}} \asymp n^2 q_n^{-\tilde{s}}/\log p_n = M^{-\tilde{s}/s}(\log n)/(\log p_n).$$

Thus the conditions of statement (*i*) hold, and $\hat{\delta}_n^Q(s)$ is a boundary rate for this test. □

The first statement of Corollary 3.1.2 is useful to deduce the $\delta_n \to 0$ such that (3.2) holds for a specific test, as is demonstrated in the proof of the second statement. The second statement is particularly important in that it establishes the concept of adaptive optimality among tests based on tapering, and characterizes the associated performance bound *via* the sequence $\hat{\delta}_n^Q(s)$. Observe that the optimal adaptive rate identified in that statement is slower than Spokoiny's adaptive-minimax rate, $\hat{\delta}_n^{AM}(s) = o(\hat{\delta}_n^Q(s))$. This of course means there are tests that would asymptotically outperform any test based on tapering in the adaptive context, such as Spokoiny's (1996) test based on thresholding or Fan's (1996) adaptive Neyman test. Nevertheless, within the class of tests based on tapering, the test based on $Q_n^{OPT}$, with $p_n$ as in Corollary 3.1.2.*ii*, is adaptively optimal. Moreover, although the condition $\log p_n \asymp \log n$ does not allow $p_n \to \infty$ at an arbitrarily fast rate, it nevertheless gives the dimensionality parameter fairly wide leeway. This property makes such a test particularly suited for use in FDA.

Let us also remark that adaptive-optimality can be established in the same manner as in the proof of Corollary 3.1.2.*ii* for the class of test statistics (1.4) with $w_{n,j} = \{j(\log j)^\gamma\}^{-1/2}$ such that $\gamma < 1$. The setting $\gamma = 0$, the special case that defines $Q_n^{OPT}$, is preferred for simplicity.

## 4. Comparisons among high-dimensional tests

In this section, the various high-dimensional tests alluded to in this article are described in detail and evaluated by simulation, with the goal of comparing the power-properties of tapering, truncation, and thresholding mechanisms.

The tests are defined here in the context of the discrete model (1.2), taking $\nu = 1$, and are compared assuming independent, Gaussian errors with known variances. Indeed, most of the available theoretical results are derived in this context, including those of Ingster (1993) and Spokoiny (1996). (Yet Fan and Lin, 1998, and Fan and Huang, 2001, establish robustness of the good power-properties of Fan's, 1996, adaptive Neyman test under weaker assumptions.) To reflect the setting $\nu = 1$, this section will revise notation to rewrite $\boldsymbol{Y}_{n,j} = Y_{n,j}$, $\boldsymbol{\theta}_{n,j} = \theta_{n,j}$, and $\boldsymbol{e}_{n,j} = e_{n,j}$ in (1.2). Each test is defined below by stating a test statistic; it should be assumed the corresponding test itself rejects the null hypothesis when the test statistic exceeds a fixed cutoff.

### 4.1. Tests based on tapering

A test based on tapering that is known to achieve Ingster's minimax rate is defined by the test statistic $FZZ_n = n \sum_{j=1}^{p_n} w_{n,j} Y_{n,j}^2$, with weights given by $w_{n,j} = 1 - j^{4s}\xi_n^2/(1 + j^{2s}\xi_n)^2$ and $\xi_n = n^{-4s/\tilde{s}}$. Optimal performance of this test was first deduced in Fan, Zhang, and Zhang (2001) for a slight variation in which the test statistic is expressed as an infinite quadratic form, $FZZ_n^\infty = n \sum_{j=1}^{\infty} w_{n,j} Y_{n,j}^2$. However, Spitzner (2008A, ex. 2) clarifies that such performance is retained for the finite-sum version, provided $n^{2/\tilde{s}} = O(p_n)$.

A minimax test based on a more simple quadratic form is studied in Ingster (1993). It is defined by the unweighted statistic $UWQ_n = n \sum_{j=1}^{p_n} Y_{n,j}^2$, and has $\hat{\delta}_n^M(s) = n^{-2s/\tilde{s}}$ as a boundary rate, provided $p_n \asymp n^{2/\tilde{s}}$. A pitfall of this test is that it is extremely sensitive to the rate at which $p_n$ increases, which puts it at a serious disadvantage in the FDA context. For instance, if $p_n \asymp n^{2/\gamma}$ for $\gamma \neq \tilde{s}$ it is possible to find sequences $\delta_n = n^{-(2-t)s/\tilde{s}}$ with $t > 0$ that fail to satisfy the rates-of-testing criterion (3.2). This is a drastic degradation in performance, well beyond that incurred to achieve adaptive-minimaxity or optimal adaptive performance among tests based on tapering. The power of the test based on $UWQ_n$ is invariant on the contours of $\sum_{j=1}^{p_n} \theta_{n,j}^2$, a property that is sometimes touted as a practical advantage; for present purposes this property will be exploited to calibrate the simulation design.

Note that tests based on $FZZ_n$, $FZZ_n^\infty$, and $UWQ_n$ each require a precise specification of the parameter *s* to achieve minimax performance. Consequently, since *s* can rarely be specified exactly, each is difficult to implement in practice. In the present exercise, these particular tests should be regarded not as candidates for practical usage, but as conceptual references against which to compare the properties of other tests.

Two other tests based on tapering will also be evaluated. One is, of course, that based on $Q_n^{OPT}$, an adaptively optimal test of this class. The other is defined by the test statistic $CVM_n = n \sum_{j=1}^{p_n} j^{-2} Y_{n,j}^2$, which is motivated by the classical Cramér-von Mises statistic of goodness-of-fit testing. The Cramér-von Mises statistic is traditionally expressed as in integrated squared-difference between empirical and hypothesized distribution functions, but it is studied in Eubank and LaRiccia (1992) *via* the representation $n \sum_{j=1}^{p_n/2} j^{-2}(Y_{n,2j-1}^2 + Y_{n,2j}^2)$, giving rise to the statistic $CVM_n$ studied here. This test is included as an example of a widely-used high-dimensional test that exhibits rather poor performance in an FDA context. For instance, it can be shown that, due to its strong down-weighting of higher-indexed components, its best possible boundary rates are $\hat{\delta}_n = n^{-(2-t)s/\tilde{s}}$ for $t = 2\{1 - \tilde{s}/(\tilde{s}+3)\} > 0$.

## 4.2. Adaptive-minimax tests based on truncation and thresholding

Two tests that achieve Spokoiny's adaptive-minimax rate shall be considered, one based on truncation and the other on thresholding. The test based on truncation is Fan's (1996) adaptive Neyman test, whose test statistic is $AN_n = \max_{k=1,\ldots,p_n}(N_{n,k} - k)/\sqrt{k}$, for which $N_{n,k} = n \sum_{j=1}^{k} Y_{n,j}^2$. Sometimes this test is interpreted *via* the scheme known as "Neyman's truncation," which describes the test statistic by $N_{n,\hat{k}_n}$, viewing $\hat{k}_n$ as a data-driven diagnostic; to yield $AN_n$, $\hat{k}_n$ is that $k$ which maximizes the standardized sums $(N_{n,k}-k)/\sqrt{k}$. Other choices of $\hat{k}_n$ are discussed in Raynor and Best (1989), Eubank and Hart (1992), Eubank and LaRiccia (1992), Inglot and Ledwina (1996), Eubank (2000), Aerts, Claeskens, and Hart (2000), Claeskens and Hjort (2004), and elsewhere. Among them, the test based on $AN_n$ has the most well established rates-of-testing properties and compares the most favorably in the empirical-power investigation of Spitzner (2006), which is similar to the investigation carried out in what follows. It is shown in Fan, Zhang, and Zhang (2001) to achieve Spokoiny's adaptive-minimax rate, at least for the case where $p_n = n$.

The test studied here based on thresholding is introduced in Spokoiny (1996). This test is readily applicable under the discrete model (1.2), but its formulation takes the continuous model (1.1) to have been translated to (1.2) using wavelet decomposition rather than Fourier decomposition. For present purposes it is unnecessary to describe the complete details of that translation, except to note that its organization of component-subscripts uses "wavelet indexing:" $j = j(k,l)$ denotes the *l*'th component at the *k*'th level of resolution: $k = 0, 1, \ldots$ and $l = 1, \ldots, 2^k$; $j(0,1) = 1$, $j(k+1, l) = j(k,l) + 2^k$, and $j(k, l+1) = j(k,l) + 1$.

The test statistic is constructed in two stages. In the first stage, a class of tests indexed by $s > 1/2$ is constructed as follows. Define

$$HT_n(s) = \sum_{k=0}^{k_{*n}(s)} \sum_{l=1}^{2^k} \left[ nY_{n,j(k,l)}^2 - 1 \right] \qquad (4.1)$$
$$+ \sum_{k=k_{*n}(s)+1}^{k_n^*} \sum_{l=1}^{2^k} \left[ nY_{n,j(k,l)}^2 I\{\sqrt{n}Y_{n,j(k,l)} > \xi_{n,k}(s)\} - \mu_{HT}(\xi_{n,k}(s)) \right]$$

for parameters $k_{*n}(s)$, $k_n^*$, and a "hard-thresholding" parameter $\xi_{n,k}(s)$, where $\mu_{HT}(\xi) = E[\eta^2 I\{|\eta| > \xi\}]$ with $\eta$ following a standard-normal distribution. The settings applied here take $k_n^* = \lceil \log_2 n \rceil$, $k_{*n}(s) = \lceil (2s+1/2)^{-1} \log_2 n \rceil$, and $\xi_{n,k}(s) = \sqrt{(k - k_{*n}(s) + 8) \log 2}$, where $\lceil x \rceil$ denotes the smallest integer still to exceed $x$. In the second stage, a range $s_* < s < s^*$ is assumed to have been specified, over which adaptively-optimal performance is to be achieved. The $HT_n(s)$ corresponding to $s$ in that range are combined into the test statistic is $\overline{HT}_n = \max_{s_* < s < s^*} HT_n(s)/c_{n,\alpha}(s)$, where $c_{n,\alpha}(s)$ is the cutoff for a size-$\alpha$ test based on $HT_n(s)$. (Note there is only a finite number of distinct $HT_n(s)$ across $s_* < s < s^*$.)

The settings for $k_{*n}(s)$ and $\xi_{n,k}(s)$ given above are not those originally specified in Spokoiny (1996), but are modified versions suitable for the non-asymptotic context considered here. The original settings are $k_{*n}(s) = \lceil (2s + 1/2)^{-1} \log_2(Mn) \rceil$, where $M$ is as in $\mathcal{B}_{s,M}$, and the hard-thresholding parameter is $\xi_{n,k}(s) = 4\sqrt{(k - k_{*n}(s) + 8) \log 2}$. With these, Spokoiny (1996) establishes that each test based on $HT_n(s)$ achieves Ingster's minimax rate, and the test based on $\overline{HT}_n$ achieves Spokoiny's adaptive-minimax rate. However, Abramovich et al. (2002) remark that in the original setting for $\xi_{n,k}(s)$ the leading constant of four is "unreasonably high" for "finite sample situations," presumably referring to data configurations one would tend to work with in practice, and proposed it be replaced with one, as is done here. The original setting for $k_{*n}(s)$ is problematic in that it depends on $M$. Abramovich et al. (2002) handle this with application-specific, ad hoc adjustments; the setting here matches the ratio $k_n^*/k_{*n}(s)$ to the limit of what it would be under the original setting.

### 4.3. An empirical comparison

The power of the tests described in Sections 4.1 and 4.2 are now compared in a simulation exercise involving two classes of high-dimensional alternatives. The first class consists of "spiked" alternatives, for which $\theta_j$ is taken to have $\theta_j = \sqrt{\lambda}$ for $j = j_0$ and $\theta_j = 0$ otherwise, where $j_0$ is a subscript that indexes the class, and $\lambda$ is set to calibrate each alternative so that the power of the test based on $UWQ_n$ is 0.4. The second class consists of "smooth" alternatives, which are parameterized according to $\theta_j = \sqrt{\lambda}\{1 - j/(p_n + 1)\}^d/c$ where $c^2 = \sum_{j=1}^{p_n} \theta_j^2$, $d = \{\log 0.2/\log(1-b) - 1\}/2$ for $0 < b < 1$, and $\lambda$ is as in the spiked class of alternatives. This formula is derived from the inverse-CDF of the beta distributions, and is such that the partial sum $\sum_{j=1}^{J} \theta_j^2$ is approximately 80% of its value at $J = p_n$ when $J$ is approximately $100b\%$ of $p_n$. The class is indexed by the parameter $b$. Examples of each type of alternative are displayed as upper-right insets in the panels of Figure 4.

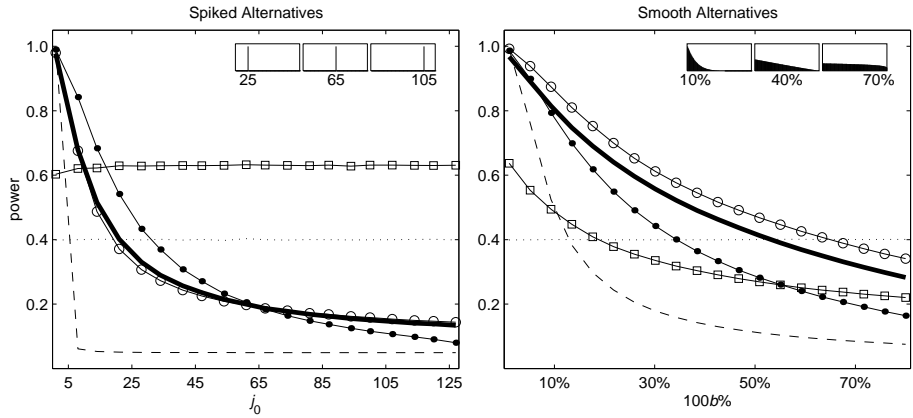

Fig 4. *Simulated power for high-dimensional tests against spiked and smooth alternatives. Power curves associated with $FZZ_n$ are marked by solid lines with circles (○), those with $UWQ_n$ by dotted lines, those with $Q_n^{OPT}$ by thick solid lines, those with $CVM_n$ by dashed lines, those with $AN_n$ by solid lines with dots (●), those with $\overline{HT}_n$ by solid lines with squares (□). Insets display the shapes of example alternatives, with corresponding index values.*

Each spiked-alternative satisfies the technical property of maximally compressing non-zero components into higher indices among alternatives with constant $\sum_{j=1}^{\infty} j^{2s}\theta_j^2$. Mathematical proofs of minimaxity single out these alternatives as yielding the minimum power, which is to be maximized (*cf.* th. 4 of Fan, Zhang, and Zhang, 2001, and th. 1 of Spitzner, 2008A). In other words, the class of spiked alternatives represent those alternatives that are the hardest to distinguish. These may be of primary interest in some specialized FDA applications, such as those involving PET-fMRI images (*cf.* Abramovich *et al.*, 2002). On the other hand, alternatives in the smooth class are idealized representations of those of primary concern in more typical FDA testing problems. For instance, in the Canadian temperature example, interest centers on such large-scale attributes as the average yearly temperature and the differential between winter and summer temperatures, attributes that are expressed in the lower-indexed $\theta_j$. The smooth alternatives reflect this and similar situations by expressing the departure from the null hypothesis mainly through these low-indexed components.

The simulation design is such that the models examined each have dimensionality $p_n = 127$, which was selected to accommodate seven levels of resolution in the wavelet-indexing scheme, at $k = 0, 1, \ldots, 6$. (This value is not atypical in FDA applications; *e.g.*, $p_n = 100$ in Fan and Lin, 1998, and $p_n = 124$ in Spitzner, Marron, and Essick, 1998.) The sample-size parameter is set to $n = 64$, which was chosen since the statistic $\overline{HT}_n$ then has $k_n^* = 6$, so that the test is sensitive to all available levels of resolution.

Simulated power is calculated at twenty alternatives from each of the spiked- and smooth-alternative classes. The specific alternatives evaluated are such that their corresponding class indices are roughly evenly-spaced across the ranges $j_0 = 1, \ldots, p_n$ for the spiked class and $0.01 \leq b \leq 0.80$ for the smooth class, including the stated endpoints.

Some intricacy is required to deduce reasonable specifications for the parameters of the Sobolev geometry, $\mathcal{B}_{s,M}$, which are needed to construct the statistics $FZZ_n$ and $\overline{HT}_n$. For spiked alternatives, the bound in $\mathcal{B}_{s,M}$ is set to $M = \lambda(p_n+1)$, so that the Sobolev norm of the alternative indexed at $j_0 = p_n+1$ would have $\sum_{j=1}^\infty j^{2s}\theta_j^2 = M$ at $s = 1/2$, if such a setting for $s$ were allowed. Then fixing $M$ at this value, the parameter $s$ is set at different values for different alternatives: for $j_0 = 2, \ldots, p_n$, the parameter is $s = \log(p_n + 1)/\{2\log j_0\}$, which solves $\sum_{j=1}^\infty j^{2s}\theta_j^2 = M$; for $j_0 = 1$, $s$ is set to its value at $j_0 = 2$. For smooth alternatives, $M$ is first set to the value $\sum_{j=1}^{p_n} j^{2s}\theta_j^2$ obtained with $s = 1/2$ and the $\theta_j$ defined by indexing at $b = 0.81$. The parameter $s$ is then set numerically at each alternative evaluated, for which $b \leq 0.80$, to the value that solves $\sum_{j=1}^\infty j^{2s}\theta_j^2 = M$.

The bounds $s_*$ and $s^*$, which are required to construct $\overline{HT}_n$, are taken as the lower and upper values of $s$ calculated in the scheme above, treating the two classes of alternatives separately. The selected spiked alternatives yield the bounds $s_* = 0.5008$ (at $j_0 = 127$) and $s^* = 1.1667$ (at $j_0 = 1$); the selected smooth alternatives yield $s_* = 0.5017$ (at $b = 0.80$) and $s^* = 2.4680$ (at $b = 0.01$). Under these settings, the parameter $k_{*n}(s)$ varies within a very narrow range: $k_{*n}(s) = 4, 5$ for spiked alternatives, and $k_{*n}(s) = 3, 4, 5$ for smooth alternatives. Consequently, the simulated power of the test based on $\overline{HT}_n$ is nearly identical to that of each test based on $HT_n(s)$.

Simulated power curves are displayed in the panels of Figure 4 for tests based on $FZZ_n$, $UWQ_n$, $Q_n^{OPT}$, $CVM_n$, $AN_n$, and $\overline{HT}_n$. Test-statistic cutoffs and the quantities $\mu_{HT}(\xi_{n,k}(s))$ of the statistic $\overline{HT}_n$ were calculated by simulation, and all tests were carried out at the $\alpha = 0.05$ level. Every simulation used a minimum of 250,000 iterations. The left panel of Figure 4 displays results for the spiked alternatives, and the right panel those of the smoothed alternatives. One should find that the test based on $UWQ_n$ serves well as a basis of comparison for the other tests, observing that its simulated power curve is near-constant at the value 0.4.

Examining first the results for spiked alternatives displayed in the left panel, the benefit of the test based on $\overline{HT}_n$ against spiked alternatives is glaring. Whereas the simulated power curve of every other test follows the same pattern of starting out high at low values of $j_0$ then dropping sharply and later evening out well below 0.4, that of $\overline{HT}_n$ is nearly flat at 0.6, representing a consistent 50% increase in power above that of the test based on $UWQ_n$. Among the remaining tests, the simulated power curves of the tests based on $FZZ_n$, $Q_n^{OPT}$, and $AN_n$ appear quite similar, and reflect performance that is not altogether poor, especially at the lower values of $j_0$. The test based on $CVM_n$, however, exhibits exceptionally poor performance.

The picture changes drastically for the smooth alternatives. In the right panel, the considerable advantage of the test based on $\overline{HT}_n$ against spiked alternatives disappears; its simulated power, in fact, appears generally low relative to the other tests, and the same can be said of the test based on $AN_n$, although to a lesser degree. To put it another way, it is the tests based on *tapering* that exhibit superior performance in this context, specifically those based on $FZZ_n$ and $Q_n^{OPT}$. The test based on $FZZ_n$ outperforms that based on $Q_n^{OPT}$, and this is not altogether surprising given that the former satisfies a stronger asymptotic property. Simulated power of the test based on $CVM_n$ is still quite poor, although not as dismal as it is for spiked alternatives.

Though the tests based on truncation and thresholding are known to have superior asymptotic performance, these empirical results suggest quite clearly that the higher-order asymptotic factors wash out slowly, and may very well have a substantial influence in data configurations one is likely to encounter in practice. All of the tests are shown to be sensitive to the *shape* of the alternative, and the test based on thresholding appears especially sensitive in this regard. These observations substantiate a recommendation to the analyst that when choosing among high-dimensional test in FDA it is prudent to consider the types of departures from the null hypothesis one is most interested in detecting. The test based on $Q_n^{OPT}$ may best suit the goals of the analysis, despite its asymptotic inferiority to truncation and thresholding, and it likely will in typical FDA applications where the primary interest is in the large-scale attributes of the model.

Let us briefly return to an observation made in the Canadian temperature example. As has been suggested, the focus on large-scale attributes puts the goals of the analysis in line with detecting smooth-shaped alternatives, and, in light of the current simulation results, the test based on $F_{global}$ is best suited to that purpose. However, recall that the test rejected H$_0$ on the example data, and exploratory analysis indicated the presence of a spiked-shaped alternative. Thus we have a practical illustration of a test based on tapering evidently detecting an alternative among those hardest to distinguish, thereby further demonstrating the *global* power of the test.

## 5. Concluding discussion

Through an application example and both theoretical and empirical power investigations, this article has demonstrated the benefits of the tapering mechanism in FDA, and the test based on $Q_n^{OPT}$ in particular. It has been shown how tests based on tapering may be constructed on a functional linear model, using Fourier decomposition to first translate to a high-dimensional discrete model. Intuition for the discrete model and tests based on tapering have been discussed through an example analysis of the Canadian temperature data. Applying criteria from rates-of-testing theory, it has been shown that the weight setting $w_{n,j} = j^{-1/2}$ in (1.4), which defines $Q_n^{OPT}$, represents an adaptively optimal configuration among tests based on tapering. Moreover, the test based on $Q_n^{OPT}$ retains optimality over a wide range of rates at which $p_n \to \infty$, a property that is particularly advantageous in the FDA context. The discussion reiterates that any test based on tapering can be outperformed asymptotically by tests based on truncation or thresholding. Nevertheless, an empirical investigation in a non-asymptotic context has demonstrated that high-order effects may be non-negligible in practice, and that the test based on $Q_n^{OPT}$ has superior power against a class of alternatives that would be of natural interest in many FDA applications.

While this article argues in support of the tapering mechanism on the basis of its power properties, there are a number of reasons completely separate from power why an analyst would, at the outset, want to restrict himself or herself to tests based on tapering. The first is that the test statistic (1.4) may arise through a formal Bayesian construction, as a monotone transformation of a posterior null probability. This is shown in Spitzner (2008B), where the rates-of-testing framework is developed entirely within a Bayesian context. It should be noted that there are existing Bayesian constructions of the thresholding mechanism for use in *estimation*, *e.g.*, in Abramovich *et al.* (2007), which suggest the statistic $\overline{HT}_n$ might be viewed as a summary metric applied a Bayesian estimator. This falls short of producing a formal Bayesian test based on $\overline{HT}_n$, however, since it fails to represent $\overline{HT}_n$ as a monotone transformation of a posterior null probability.

Another reason that the tapering mechanism is attractive is because of its straightforward intuition. The non-expert, once appreciating the goal to exploit smoothness but still retain global power, would more quickly embrace the intuition underlying the use of tapering before truncation or thresholding. Tapering is easy to understand as it directly incorporates all components of the model, while explicitly down-weighting the high-indexed ones. The latter mechanisms are more complicated, harder to connect to the goal, and easy to apply incorrectly. For these reasons, tapering may be the preferred recommendation in consulting or interdisciplinary situations, and may offer the best insurance for continued correct implementation once the statistician's involvement in a project wanes.

Finally, the test statistic (1.4) is convenient in that it may be treated through well-known analytical approximations (see, *e.g.*, Mathai and Provost, 1992, sec. 4.6), rather than simulation. This may be a trivial advantage in straightforward testing problems, given modern computing power, but it can help enormously when developing more complicated high-dimensional statistical procedures. For instance, Spitzner and Marshall (2009) make critical use of distributional approximations of quadratic forms to develop a high-dimensional sequential monitoring procedure based on tapering.

The results of this article encourage the use of the test based on $Q_n^{OPT}$ in FDA, and establish a role for the tapering mechanism in rates-of-testing theory. It is hoped that the ideas presented here are also found helpful in clarifying the role of smoothness constraints in high-dimensional problems, and in providing an accessible methodology with which to exploit them.

## Acknowledgements

The author is grateful for the comments of the editor, associate editor, and two referees, whose suggestions substantially improved the presentation.

## References

Abramovich, F., Antoniadis, A., Sapatinas, T., and Vidakovic, B. (2002). Optimal testing in functional analysis of variance models. *Technical Report* G02/2002. http://zoe.bme.gatech.edu/∼bv20/dpgatech.html. (See also, *International Journal of Wavelets, Multiresolution, and Information Processing*, 2:323–349.) MR2104869

Abramovich, F., Grinshtein, V., and Pensky, M. (2007). On optimality of Bayesian testimation in the normal means problem. *Annals of Statistics* **35** 2261-2286. MR2363971

Adams, R. A., and Fournier, J. J. F. (2003). *Sobolev Spaces*. Academic Press, Boston. MR2424078

Aerts, M., Claeskens, G., and Hart, J. (2000). Testing lack of fit in multiple regression. *Biometrika* **87** 405-424. MR1782487

Brockwell, P. J., and Davis, R. A. (1991). *Time Series: Theory and Methods*, 2nd ed. Springer, New York. MR1093459

Claeskens, G., and Hjort, N. L. (2004). Goodness of fit via nonparametric likelihood ratios. *Scandinavian Journal of Statistics* **31** 487-513. MR2101536

Eubank, R. (2000). Testing for no effects by cosine methods. *Scandinavian Journal of Statistics* **27** 747-763. MR1804174

Eubank, R., and Hart, J. (1992). Testing goodness-of-fit in regression via order selection criteria. *Annals of Statistics* **20** 1412-1425. MR1186256

Eubank, R., and LaRiccia, V. N. (1992). Asymptotic comparison of Cramér-von Mises and nonparametric function estimation techniques for testing goodness of fit. *Annals of Statistics* **20** 2071-2086. MR1193326

Fan, J. (1996). Tests of significance based on wavelet thresholding and Neyman's truncation. *Journal of the American Statistical Association* **91** 674-688. MR1395735

Fan, J., and Huang, L.-S. (2001). Goodness-of-fit testing for parametric regression models. *Journal of the American Statistical Association* **96** 640-652. MR1946431

Fan, J., and Lin, S.-K. (1998). Test of significance when data are curves. *Journal of the American Statistical Association* **93** 1007-1021. MR1649196

Fan, J., Zhang, C., and Zhang, J. (2001). Generalized likelihood ratio statistics and Wilks phenomenon. *Annals of Statistics* **29** 153-193. MR1833962

Inglot, T., and Ledwina, T. (1996). Asymptotic optimality of data driven Neyman's tests. *Annals of Statistics* **24** 1982-2019. MR1421157

Ingster, Y. I. (1993). Asymptotically minimax hypothesis testing for nonparametric alternatives I-III. *Mathematical Methods of Statistics* **2** 85-114; **3** 171-189; **4** 249-268.

Gayraud, G., and Pouet, C. (2005). Adaptive minimax testing in the discrete regression scheme. *Probability Theory and Related Fields* **133** 531-558. MR2197113

Horowitz, J. L., and Spokoiny, V. G. (2001). An adaptive, rate-optimal test of a parametric mean-regression model against a nonparametric alternative. *Econometrica* **69** 599-631. MR1828537

Lepski, O. V., and Spokoiny, V. G. (1999). Minimax nonparametric hypothesis testing: the case of an inhomogeneous alternative. *Bernoulli* **5** 333-358. MR1681702

Mathai, A. M., and Provost, S. B. (1992). *Quadratic Forms in Random Variables: Theory and Applications.* M. Dekker, New York. MR1192786

Rayner, J. C. W., and Best, D. J. (1989). *Smooth Tests of Goodness-of-Fit.* Oxford University Press, Oxford. MR1029526

Ramsay, J. O., and Silverman, B. W. (2005). *Functional Data Analysis*, 2nd ed. Springer, New York. (Website: www.functionaldata.org) MR2168993

Seber, G. A. F. (2003). *Linear Regression Analysis*, 2nd ed. Wiley, New York. MR0436482

Spitzner, D. J. (2006), Use of goodness-of-fit procedures in high-dimensional testing. *Journal of Statistical Computation and Simulation* **76** 447-457. MR2224364

Spitzner, D. J. (2008A). Testing in functional data analysis using quadratic forms. *Technical Report Series* 08-01. University of Virginia, Department of Statistics, Charlottesville, VA.

Spitzner, D. J. (2008B). An asymptotic viewpoint on high-dimensional Bayesian testing. *Bayesian Analysis.* **3** 121-160.

Spitzner, D. J., Marron, J. S., and Essick, G. (2003). Mixed-model functional ANOVA for studying human tactile perception. *Journal of the American Statistical Association* **98** 263-272. MR1982574

Spitzner, D. J., and Marshall, J. (2009). Directed spatio-temporal monitoring of disease incidence rates. In preparation.

Spitzner, D. J., and Woodall, W. H. (2003). High-dimensional directed testing for monitoring functional profiles. *Proceedings of the American Statistical Association* 225-236.

Spokoiny, V. G. (1996). Adaptive hypothesis testing using wavelets. *Annals of Statistics* **24** 2477-2498. MR1425962